\documentclass[12pt]{article} 
\usepackage{amsfonts,latexsym,color}
\input{psfig.sty}

\newcounter{environment}[section]
\renewcommand{\theenvironment}{%
\arabic{section}.\arabic{environment}}

{\begin{rm}\refstepcounter{environment}{\textbf\theenvironment\
\bf Definition.~~}}%
{\end{rm}}

\newenvironment{example}%
{\begin{rm}\refstepcounter{environment}{\textbf\theenvironment\
\bf Example.~~}}%
{\end{rm}}

{\begin{rm}\refstepcounter{environment}{\textbf\theenvironment\
\bf Conjecture.~~}}%
{\end{rm}}

{\begin{rm}\refstepcounter{environment}{\textbf\theenvironment\
\bf Proposition.~~}}%
{\end{rm}}

{\begin{rm}\refstepcounter{environment}{\textbf\theenvironment\
\bf Proposition and Definition.~~}}%
{\end{rm}}

\newenvironment{theorem}%
{\begin{rm}\refstepcounter{environment}{\textbf\theenvironment\
\bf Theorem.~~}}%
{\end{rm}}

{\begin{rm}\refstepcounter{environment}{\textbf\theenvironment\
\bf Corollary.~~}}%
{\end{rm}}

{\begin{rm}\refstepcounter{environment}{\textbf\theenvironment\
\bf Lemma.~~}}%
{\end{rm}}

\begin{document}
\newcommand{\qbc}[2]{ {\left [{#1 \atop #2}\right ]}}
\newcommand{\anbc}[2]{{\left\langle {#1 \atop #2} \right\rangle}}
\newcommand{\be}{\begin{enumerate}} \newcommand{\ee}{\end{enumerate}}
\newcommand{\beq}{\begin{equation}} \newcommand{\eeq}{\end{equation}}
\newcommand{\bea}{\begin{eqnarray}} \newcommand{\eea}{\end{eqnarray}}
\newcommand{\beas}{\begin{eqnarray*}}
\newcommand{\eeas}{\end{eqnarray*}} 
\newcommand{\zz}{\mathbb{Z}}
\newcommand{\pp}{\mathbb{P}} 
\newcommand{\nn}{\mathbb{N}}
\newcommand{\qq}{\mathbb{Q}} 
\newcommand{\rr}{\mathbb{R}}
\newcommand{\rrn}{\rr^{1,n}}
\newcommand{\bm}[1]{{\mbox{\boldmath $#1$}}}
\newcommand{\bmu}{\bm{u}}
\newcommand{\bmv}{\bm{v}}
\newcommand{\bmp}{\bm{p}}
\newcommand{\bmx}{\bm{x}}
\newcommand{\bmy}{\bm{y}}
\newcommand{\bmq}{\bm{q}}
\newcommand{\ca}{{\cal A}}
\newcommand{\cb}{{\cal B}}
\newcommand{\ppk}{\bm{p_1},\dots,\bm{p_k}} 
\newcommand{\hz}{\hat{0}}
\newcommand{\sn}{\mathfrak{S}_n} 
\newcommand{\sk}{\mathfrak{S}_k} 
\newcommand{\cs}{{\cal S}}
\newcommand{\ds}{\displaystyle} 
\newcommand{\fs}{\mathfrak{S}}
\newcommand{\st}{\,:\,} \newcommand{\maj}{\mathrm{maj}}
\newcommand{\tx}{(t,\bm{x})}
\newcommand{\tr}{\textcolor{red}} \newcommand{\tb}{\textcolor{blue}}
\newcommand{\tg}{\textcolor{green}}
\newcommand{\tm}{\textcolor{magenta}}
\newcommand{\tbn}{\textcolor{brown}}
\newcommand{\tp}{\textcolor{purple}}
\newcommand{\tn}{\textcolor{nice}}
\newcommand{\tor}{\textcolor{orange}}

\definecolor{brown}{cmyk}{0,0,.35,.65}
\definecolor{purple}{rgb}{.5,0,.5}
\definecolor{nice}{cmyk}{0,.5,.5,0}
\definecolor{orange}{cmyk}{0,.35,.65,0}

\begin{centering}
\textcolor{red}{\Large\bf Ordering Events in Minkowski Space}\\[.2in] 
\textcolor{blue}{Richard P. Stanley}\footnote{Partially supported by
  NSF grant \#DMS-9988459.}\\ 
Department of Mathematics\\
Massachusetts Institute of Technology\\
Cambridge, MA 02139\\
\emph{e-mail:} rstan@math.mit.edu\\[.2in]
\textcolor{magenta}{version of 26 May 2005}\\[.2in]

\textbf{Abstract}\\[.1in]
\end{centering}
Let $\ppk$ be $k$ points (events) in $(n+1)$-dimensional Minkowski
space $\rrn$. Using the theory of hyperplane arrangments and chromatic
polynomials, we obtain information on
the number of different orders in which the 
events can occur in different reference frames if the events are
sufficiently generic. We consider the question of what sets of
orderings of the points are possible and show a connection with
sphere orders and the allowable sequences of Goodman and Pollack.
\vskip 10pt
\section{Introduction.} \label{sec1}
\indent Let $\rr^{1,n}$ denote $(n+1)$-dimensional Minkowski space. We
denote 
points (or \emph{events}) in $\rr^{1,n}$ as $\bm{p}=(t,\bm{x})$, where
$t\in\rr$ is the time coordinate and $\bm{x}=(x_1,\dots,x_n)\in\rr^n$
are the space coordinates. The geometry of $\rr^{1,n}$ is defined by
the \emph{Minkowski norm}
   \beas |(t,\bm{x})|^2 & = & t^2-|\bm{x}|^2\\ & = &
     t^2 -(x_1^2+\cdots+x_n^2). \eeas
An event $\tx$ is said to be \emph{timelike} if $|\tx|^2>0$,
\emph{lightlike} if $|\tx|=0$, and \emph{spacelike} if
$|\tx|^2<0$. Two events $(s,\bm{x})$ and $(t,\bm{y})$ are
\emph{timelike separated} if their difference $(s,\bm{x})-(t,\bm{y})$
is timelike, and similarly \emph{lightlike separated} and
\emph{spacelike separated}. 
Two events are timelike separated if and only if they are
\emph{causally related}; a signal (traveling slower than $c=1$, where
$c$ denotes the speed of light) can reach one event from the other. 

Now suppose that $F'$ is a second reference frame, moving with
constant velocity $\bm{v}\in\rr^n$ with respect to the original frame
$F$. By convention an observer in the frame $F'$ measures
coordinates $(t',\bm{x'})$, synchronized so that $t=t'=0$ when the two
frames coincide. Write
  $$ \bm{v} =(\tanh \rho)\bm{u}, $$
where $\bm{u}$ is a unit vector and $\tanh\rho=|\bm{v}|<1$. The
\emph{Lorentz transformation} expresses the $F'$ coordinates
$(t',\bm{x'})$ in terms of the $F$ coordinates $\tx$. All we
need here is the formula for $t'$:
  $$ t' = (\cosh \rho)t-(\sinh\rho)\bm{x}\cdot\bm{u}, $$
where $\bm{x}\cdot\bm{u}$ is the ordinary dot product of two vectors
in $\rr^n$. 

It is easy to verify from the Lorentz transformation that two timelike
separated events occur in the same order for any observers (i.e., in
any reference frame $F'$). (Otherwise, in fact, causality would be
violated.) On the other hand, two spacelike separated events can
always occur in either order in suitable reference frames. These facts
for two events suggest generalizing to more events.

\medskip
\textbf{Main problem.} Given $k$ events in $\rr^{1,n}$ in what
different orders can they occur for different observers? How many such
orders are there?
\medskip

For instance, given three events $\bm{p_1},\bm{p_2},\bm{p_3}$ in
$\rr^{1,1}$, we will see in Section~\ref{sec4} that there do not exist
three observers for which the events occur in the orders $\bm{p_1} <
\bm{p_2} < \bm{p_3}$, (i.e., $\bm{p_1}$ before $\bm{p_2}$ before
$\bm{p_3}$ in time) $\bm{p_2}<\bm{p_3}<\bm{p_1}$, and
$\bm{p_3}<\bm{p_1}<\bm{p_2}$.  On the other hand, given any two
permutations $\pi_1$ and $\pi_2$ of $[k] = \{1,2,\dots,k\}$, there
exist $k$ events $\ppk\in\rr^{1,1}$ and two observers $F_1, F_2$ such
that for $F_i$ the events occur in the order $\pi_i$, $1\leq i\leq 2$.
In general, we write ${\cal O}(\bm{p_1},\dots,\bm{p_k})$ for the
number of different orders in which the events
$\bm{p_1},\dots,\bm{p_k}$ occur for different observers.

\section{Generic spacelike separated events.} \label{sec2}
\indent In this section we assume that all the events $\bm{p_i}$ are
spacelike separated. It is easy to construct such events, e.g., choose
any events with different space coordinates and dilate their space
coordinates by a sufficiently large factor. Let $(t_i,\bm{x_i})$ be
the coordinates of $\bm{p_i}$ with respect to a fixed observer. An
observer moving at velocity $\bm{v}=(\tanh \rho)\bm{u}$ sees
$\bm{p_i}$ occur at time
  $$ t'_i = (\cosh\rho)t_i-(\sinh\rho)\bm{x_i}\cdot\bm{u}. $$
Hence $\bm{p_i}$ occurs simultaneously to $\bm{p_j}$ for this observer
if $t'_i=t'_j$, i.e., 
  $$ (\cosh\rho)t_i-(\sinh\rho)\bm{x_i}\cdot\bm{u} =
     (\cosh\rho)t_j-(\sinh\rho)\bm{x_j}\cdot\bm{u}. $$
Equivalently,
  $$ t_i-t_j=(\bm{x_i}-\bm{x_j})\cdot \bm{v}. $$
The set of such velocities $\bmv$ forms a hyperplane in $\rr^n$. The
different sides of this hyperplane determine whether $\bm{p_i}$ occurs
before or after $\bm{p_j}$. In order for the ordering of points
determined by a region $R$ to correspond to an actual observer, we must
have $|\bm{v}|<1$ for some $\bm{v}\in  R$.  Thus we obtain the
following result. 

\medskip
\begin{theorem} \label{main}
\emph{The number ${\cal O}(\bm{p_1},\dots,\bm{p_k})$
  is equal to the number of regions $R$ of
  the hyperplane arrangement $\ca = \ca(\bm{p_1},\dots,\bm{p_k})$ with
  hyperplanes $H_{ij}$ given by} 
  \beq t_i-t_j = (\bm{x_i}-\bm{x_j})\cdot \bm{v},\ \ 
     1\leq i<j\leq k, \label{eq:hyp} \eeq
\emph{such that $|\bm{v}|<1$ for some $\bm{v}\in R$.} 
\end{theorem}

For basic results about hyperplane arrangements and their number of
regions, see e.g.\ \cite{o-t}\cite{rs:hyp}. In general we denote the
number of regions of an arrangement $\ca$ by $r(\ca)$.

In general, there is no simple formula for the number of regions of
the arrangements ${\cal A}(\bm{p_1},\dots,\bm{p_k})$. There are
general formulas from the theory of arrangements for $r({\cal A})$ for
any (finite) arrangement ${\cal A}$ (e.g., equation (\ref{eq:zas})),
but such formulas do not shed much further light per se on the main
problem. If, however, we assume that $\bm{p_1},\dots, \bm{p_k}$ are
generic (in a sense to be made precise), then more can be said.
Moreover, for fixed $k$ the quantity $r(\ca(\bm{p_1},\dots,\bm{p_k}))$
is maximized when $\bm{p_1},\dots, \bm{p_k}$ are generic.

We now review a fundamental result of Zaslavsky
\cite[Thm.~2.3.21]{o-t} \cite[Thm.~2.5]{rs:hyp}\cite{zas} for
computing the number of regions of an arrangement ${\cal A}$ in
$\rr^n$ . The \emph{intersection poset} $L_\ca$ of $\ca$ is the set of
all \emph{nonempty} intersections of hyperplanes in $\ca$, ordered by
reverse inclusion. We always include the ambient space $\rr^n$ as the
bottom element $\hz$ of $L_\ca$.
The \emph{M\"obius function} $\mu:L_\ca\rightarrow \zz$ of $L_\ca$ is
defined recursively by $\mu(\hz)=1$, and for all $y>\hz$ in $L_\ca$,
  $$ \sum_{x\leq y} \mu(x) = 0. $$
(Usually $\mu$ is defined on \emph{intervals} of $L_\ca$, not
elements, so our $\mu(x)$ corresponds to $\mu(\hz,x)$.) The
\emph{characteristic polynomial} $\chi_\ca(t)$ is defined by
  $$ \chi_\ca(t) = \sum_{x\in L_\ca}\mu(x)t^{\dim(x)}, $$
where $\dim(x)$ refers to the dimension of $x$ as an affine subspace
of $\rr^n$. Zaslavsky's theorem then states that
  \beq r(\ca) = (-1)^n\chi_\ca(-1). \label{eq:zas} \eeq
Since it is known that $(-1)^{n-\dim x}\mu(x)>0$, equation
(\ref{eq:zas}) can be restated as
  \beq r(\ca) = c_0+c_1+\cdots+c_n, \label{eq:csum} \eeq
where $\chi_\ca(t)=c_0t^n-c_1t^{n-1}+\cdots+(-1)^nc_n$. 

\medskip
\begin{example}
An example relevant to our results is the \emph{braid arrangement}
$\cb_k$ in $\rr^k$, with hyperplanes
  $$ z_i=z_j,\ \ 1\leq i<j\leq k. $$
The intersection poset $L_{\cb_k}$ of $\cb_k$ is just the lattice of
partitions of the set $[k]=\{1,2,\dots,k\}$, ordered by
refinement. (See \cite[Exam.~3.10.4]{ec1} for a discussion of this
lattice.) A 
partition such as 134--26--5 (i.e., the partition with blocks
$\{1,3,4\}$, $\{2,6\}$, $\{5\}$) corresponds to the intersection
$z_1=z_3=z_4$, $z_2=z_6$. The characteristic polynomial of $\cb_k$ is
given by
  \bea \chi_{\cb_k}(t) & = &  t(t-1)\cdots (t-k+1)\nonumber \\
   & = & \sum_{i=1}^k (-1)^{k-i}c(k,i)t^i, \label{eq:cki} \eea
where $c(k,i)$ denotes a \emph{signless Stirling number of the first
kind} (the number of permutations of $[k]$ with $i$ cycles)
\cite[{\S}1.3]{ec1}. We set $c(k,i)=0$ if $i<1$ or $i>k$.
\end{example}

We now define what we mean for the spacelike separated events $\ppk$
in $\rr^{1,n}$ to be generic. We state the condition in terms of the
hyperplanes $H_{ij}$ given by (\ref{eq:hyp}). These hyperplanes are
defined for $i<j$, but we extend this definition to all $i\neq j$, so
$H_{ij}=H_{ji}$. Namely, $\ppk$ are \emph{generic} if (1) no $n+1$ of
the hyperplanes $H_{ij}$ intersect (i.e., have nonempty intersection),
and (2) the minimal subsets of the $H_{ij}$'s with $m\leq n$ elements
that do intersect have the form $C=\{ H_{i_1,i_2}, H_{i_2,i_3}, \dots,
H_{i_{m-1},i_m}, H_{i_m,i_1}\}$, where $i_1,i_2,\dots,i_m$ are
distinct. Note that such sets $C$ do indeed intersect.  It is easy to
see that ``almost all'' $k$-element sequences $\ppk$ of spacelike
separated points are generic, i.e., those that aren't form a set of
measure 0 in the space of all $k$-tuples of spacelike separated points
in $\rr^{1,n}$. 

We come to the main result of this section. If $\bm{p}=(t,\bm{x})$ and
$a\in\rr$, then write $\bm{p}(a)=(t,a\bm{x})$.

\medskip
\begin{theorem} \label{thm:ssn}
\emph{Let $\ppk$ be $k$ spacelike
  separated events in $\rr^{1,n}$. Then}
  \beq {\cal O}(\ppk) \leq c(k,k)+c(k,k-1)+\cdots+c(k,k-n).
     \label{eq:rappk}  \eeq
\emph{Moreover, we have}
  \beq {\cal O}(\bm{p_1}(a),\dots,\bm{p_k}(a)) =
         c(k,k)+c(k,k-1)+\cdots+c(k,k-n) \label{eq:equal} \eeq
\emph{for sufficiently large $a$ if and only if $\ppk$ are generic.
In particular, if $n\geq k-1$ and $\ppk$ are generic, then}
   $$ {\cal O}(\bm{p_1}(a),\dots,\bm{p_k}(a)) =k! $$ 
\emph{for sufficiently large $a$, i.e., the events can occur in any
order}.  
\end{theorem}

\textbf{Proof.} Suppose that $\ppk$ are generic. Consider the
intersection poset $L_\ca$ of the 
arrangement $\ca=\ca(\ppk)$. By genericity, no $n+1$ of the
hyperplanes intersect, and the minimal subsets that do intersect have
the form $\{ H_{i_1,i_2}, H_{i_2,i_3}, \dots, H_{i_{m-1},i_m},
H_{i_m,i_1}\}$, where $m\leq n$. This is exactly the same as for the
braid arrangement $\cb_k$, where $H_{ij}$ corresponds to the
hyperplane $z_i=z_j$, except there is no condition that $m\leq n$ for
$\cb_k$. It follows that the elements of $L_\ca$ correspond to
partitions of $[k]$ with at least $k-n$ blocks, with the partial
ordering on $L_\ca$ corresponding to refinement of partitions. Hence
$L_\ca$ is isomorphic to the rank $n$ \emph{truncation} of $\Pi_k$,
i.e., the subposet of $\Pi_k$ consisting of all elements of rank at
most $n$ (where the rank of a partition $\pi$ with $j$ blocks is
$k-j$). By (\ref{eq:cki}) we have
  $$ \chi_\ca(t) = \sum_{i=0}^n (-1)^{n-i}c(k,k-i)t^{n-i}. $$
Hence it follows from Theorem~\ref{main} and (\ref{eq:csum}) that 
  $$ {\cal O}(\ppk)\leq r(\ca) = c(k,k)+c(k,k-1)+\cdots+c(k,k-n). $$
Now replacing each point $\bm{p_i}$ with $\bm{p_i}(a)$ for $a>0$
dilates the arrangement $\ca(\ppk)$ by a factor of $1/a$ and maintains
genericity. Hence for $a$ sufficiently large every region of the
dilated arrangement $\ca(\bm{p_1}(a),\dots, \bm{p_k}(a))$ intersects
the open unit ball $|\bm{v}|<1$, showing that equation
(\ref{eq:equal}) holds for $a$ sufficiently large.

It remains to show that equality cannot hold in (\ref{eq:rappk}) if
$\ppk$ are not generic. By Theorem~\ref{main} it suffices to show that
in this case, 
  $$ r(\ca(\ppk)) < c(k,k)+c(k,k-1)+\cdots+c(k,k-n). $$
Assume that $\ppk$ are any $k$ spacelike separated events in
$\rr^{1,n}$. Let $\kappa\ca=\kappa\ca(\ppk)$ denote the \emph{cone}
over $\ca$ \cite[{\S}1.2]{o-t}\cite[p.~7]{rs:hyp}, i.e., introduce a
new coordinate $u$ and define the hyperplanes of $\kappa\ca$ by
     \beas (t_i-t_j)u & = & (\bm{x_i}-\bm{x_j})\cdot \bm{v},\ \ 
     1\leq i<j\leq k\\ u & = & 0. \eeas
It is not hard to see \cite[end of {\S}4.2]{rs:hyp} that $r(\kappa\ca)
= 2r(\ca)$. Let $L$ be a linear ordering of the hyperplanes of
$\kappa\ca$. A \emph{circuit} of $\kappa\ca$ is a minimal set of
linearly dependent hyperplanes. (A set of hyperplanes is
defined to be linearly independent if their normals are linearly
independent.) A \emph{broken circuit} is a circuit with its largest
hyperplane (in the order $L$) deleted. It is an immediate consequence
of the broken circuit theorem
\cite[(6.73)]{b-o}\cite[Thm.~4.12]{rs:hyp} and equation
(\ref{eq:csum}) that  
  \beq r(\kappa\ca)=\#\{ S\subseteq\kappa\ca\st S\ \mbox{contains no
     broken circuit}\}. \label{eq:bc} \eeq
Now suppose that $\ppk$ are any $k$ spacelike separated points in
$\rrn$. Let $\bm{p'_1},\dots,\bm{p'_k}$ be $k$ generic spacelike
separated points in $\rrn$. Denote the hyperplanes of $\ca=\ca(\ppk)$
by $H_{ij}$ as in Theorem~\ref{main}, and the corresponding
hyperplanes in $\ca'=\ca(\bm{p'_1},\dots,\bm{p'_k})$ by $H'_{ij}$. Let
$J$ denote the hyperplane $u=0$ of $\kappa\ca$ and $J'$ the
corresponding hyperplane of $\kappa\ca'$. If $S\subseteq \kappa\ca$,
then let $S'=\{H'\st H\in S\}$.

Now let $C'$ be a circuit of $\kappa\ca'$. Then $C$ is a linearly
dependent subset of $\kappa\ca$. Hence if $B$ is a broken circuit of
$\kappa\ca$, then $B'$ is contained in a broken circuit of
$\kappa\ca'$. It follows from Theorem~\ref{thm:ssn} and equation
(\ref{eq:bc}) that 
  \beas 2r(\ca) & = & r(\kappa\ca)\\ & \leq & 
     r(\kappa\ca')\\ & = & 2r(\ca')\\ & = &
     c(k,k)+c(k,k-1)+\cdots+c(k,k-n). \eeas
\indent It remains to show that if $\ppk$ are not generic, then the
above inequality is strict. This is equivalent to showing that there
exists a broken circuit $B$ of $\kappa\ca$ such that $B'$ contains no
broken circuit of $\kappa\ca'$. We are free to choose any linear
ordering $L$ of $\kappa\ca$ that is convenient, and the corresponding
linear ordering $L'$ of $\kappa\ca'$ (i.e., if $H<K$ in $L$, then
$H'<K'$ in $L'$). Let $C$ be a circuit of $\kappa\ca$ such that $C'$
is not a circuit of $\kappa\ca'$. Such a circuit $C$ exists since
$\ppk$ are not generic. Thus $C'$ is linearly independent. Let 
  $$ X' =\{H'\in\kappa\ca'\st C'\cup\{H'\}\ \mbox{contains a circuit}
     \}. $$
Choose $L'$ so that all elements of $X'$ come before all elements of
$C'$. Then $C'$ contains no broken circuit of $\kappa\ca'$ with
respect to $L'$. Hence if $D$ is a broken circuit of $\kappa\ca$
contained in $C$ (e.g., we can always take $D$ to be $C$ minus its
largest element), then $D'$ contains no broken circuit of
$\kappa\ca'$, completing the proof. $\ \Box$
    
\textsc{Note.} The above argument extends to any matroid and shows the
following. (For matroid theory terminology, see e.g.\ \cite{mat}.) Let
$M$ be a (finite) matroid with characteristic polynomial
$t^m-a_1t^{m-1}+\cdots+(-1)^ma_m$. Let $N$ be a weak map image of $M$
with characteristic polynomial $t^n-b_1t^{n-1}
+\cdots+(-1)^nb_n$. Then $b_i\leq a_i$ for all $i$. This result is
essentially known \cite[Props.~7.3, 7.4]{lucas}, but since it is only
given in the case rank$(M)=\mathrm{rank}(N)$ (i.e., $m=n$) we have
provided the above proof.

\section{Timelike separated events.} \label{sec3}
\indent We consider in this section the
case where some of the events are timelike separated. Recall that if
$(s,\bm{x})$ and $(t,\bm{y})$ are timelike separated, then they occur
in the same order in all reference frames. In that case, solutions 
$\bmv$ to $s-t=(\bm{x}-\bm{y})\cdot \bmv$ satisfy $|v|>1$, so these
hyperplanes can be ignored since they are physically meaningless. Thus
let $\ppk$ be \emph{any} $k$ events in $\rr^{1,n}$. For convenience
assume no two are lightlike separated. Define the \emph{separation
  graph} $G=G(\ppk)$ to be the (undirected) graph on the vertex set
$V(G)=[k]$ with edge set
  $$ E(G)=\{ ij\st \bm{p_i},\bm{p_j}\ \mbox{are spacelike
         separated}\}. $$
The following extension of Theorem~\ref{main} is clear.

\medskip
\begin{theorem} \label{thm:ag}
\emph{Let $\ppk$ be events in $\rrn$ with separation graph $G$. Then
 ${\cal O}(\bm{p_1},\dots,\bm{p_k})$
  is equal to the number of regions $R$ of
  the hyperplane arrangement $\ca = \ca(\bm{p_1},\dots,\bm{p_k})$ with
  hyperplanes $H_{ij}$ given by} 
  $$ t_i-t_j = (\bm{x_i}-\bm{x_j})\cdot \bm{v},\ \ 
     ij\in E(G), $$
 \emph{such that $|\bm{v}|<1$ for some $\bm{v}\in R$.} 
\end{theorem}
\medskip

Let $G$ be a graph with $V(G)=[k]$. The \emph{graphical arrangement}
$\cb_G$ is the hyperplane arrangement in $\rr^k$ with hyperplanes
$z_i=z_j$ for $ij\in E(G)$. For instance, if $G=K_k$ (the complete
graph on $[k]$), then $\cb_{K_k}=\cb_k$, the braid arrangement. Let
$\chi_G(t)$ denote the \emph{chromatic polynomial} of $G$, i.e., for
$m\in\pp$, $\chi_G(m)$ is the number of ways to color the vertices of
$G$ from a set of $m$ colors such that adjacent vertices have
different colors. It is well-known
\cite[Thm.~2.4.19]{o-t}\cite[Thm.~2.7]{rs:hyp} 
that the characteristic polynomial of $\cb_G$ is given by
  $$ \chi_{\cb_G}(t) = \chi_G(t). $$
\indent Fix a graph $G$ on the vertex set $[k]$, and assume that
$\ppk$ satisfy $G=G(\ppk)$ but are otherwise generic. Then just as
for the case $G=K_k$, we have that
$L$$_{\ca(\mbox{\small{\boldmath $\ppk$}})}$ is the rank $n$
truncation of $L_{\cb_G}$. We obtain just as for Theorem~\ref{thm:ssn}
(the special case $G=K_k$) the following result. 

\medskip
\begin{theorem} \label{thm:rag}
\emph{Let $\bm{p_1},\dots,\bm{p_k}\in\rr^{1,n}$ be $k$ events in
$\rrn$ with separation graph $G$.  Let $\chi_G(t)=t^k-a_1
t^{k-1}+\cdots+(-1)^{k-1}a_{k-1}t$. Set $a_i=0$ if $i\geq k$.
Then} 
  \beq {\cal O}(\ppk)\leq 1+a_1+a_2+\cdots+a_n, \label{eq:chrobd}
  \eeq  
\emph{with equality only if $\ppk$ are generic (with respect to 
having separation graph $G$).}
\end{theorem}
\medskip

Unlike Theorem~\ref{thm:ssn} we don't necessarily have equality
holding in (\ref{eq:chrobd}) for generic $\bm{p_i}(a)$
and $a$ sufficiently large, because the transformation $\bm{p}\mapsto
\bm{p}(a)$ for large $a$ will not preserve the separation
graph. Timelike separated points will become spacelike separated. The
following problem is therefore suggested.
\smallskip

\textbf{Problem.} Characterize those graphs $G$ for which there exist
events $\ppk$ in $\rrn$ with separation graph $G$ such that 
  $$ {\cal O}(\ppk) = 1+a_1+a_2+\cdots+a_n, $$
where $\chi_G(t)=t^k-a_1 t^{k-1}+\cdots+(-1)^{k-1}a_{k-1}t$.
\smallskip

Theorems~\ref{thm:ag} and \ref{thm:rag} and the above problem suggest
the problem of characterizing separation graphs of subsets of
$\rr^{1,n}$. This problem has been considered previously, and we
briefly summarize known results. Given $\bmp=(s,\bmx)\in\rr^{1,n}$,
define the (open) \emph{future light cone} $C(\bmp)$ to consist of all
points $\bmq=(t,\bmy) \in\rrn$ such that (1) $t>s$ and (2) $\bmp$ and
$\bmq$ are timelike separated. Equivalently,
  $$ C(\bmp) = \{(t,\bmy)\in\rrn\st t-s>|\bmy-\bmx|\}, $$
a half-cone with apex $\bmp$, slope $45^\circ$, and opening in the
$t$-direction. Note that if $\bmq\in C(\bmp)$ then $C(\bmq)\subset
C(\bmp)$. Define $(s,\bmx)< (t,\bmy)$ if $s< t$ and if $(s,\bmx)$
and $(t,\bmy)$ are timelike separated. It follows that the reflexive
closure of the relation $<$ (i,e., define $\bmp\leq\bmq$ if
$\bmp<\bmq$ or $\bmp=\bmq$) is a partial order $P_n$. Any induced
subposet of $P_n$ is called a \emph{timelike poset} or \emph{causal
  poset}. Thus if $G$ is the separation graph of a finite subset $S$
of $\rrn$, then $G$ is the \emph{incomparability graph} of the
restriction of $P_n$ to $S$. In other words, the vertices of $G$ are
the elements of $S$, with an edge between two vertices if they are
incomparable in $P_n$.

It is a strong restriction on separation graphs to be incomparability
graphs. See e.g.\ \cite[{\S}3.2]{trotter} for some characterizations of
incomparability graphs. We may further ask what other conditions are
satisfied by separation graphs. Suppose $C(\bmp)$ and $C(\bmq)$ are
future light cones. Intersect them with a hyperplane $t=t_0$ for $t_0$
large. The intersections are just balls $B(\bmp)$ and $B(\bmq)$.
Moreover, $C(\bmp)\subset C(\bmq)$ if and only if $B(\bmp)\subset
B(\bmq)$. It follows that a finite poset $P$ is a timelike poset in
$\rrn$ if and only if it is an $n$-dimensional \emph{sphere order},
i.e., isomorphic to a set of spheres in $\rr^n$, ordered by inclusion
of their interiors. In fact, the concept of sphere orders originally
arose in the above context of special relativity \cite{meyer}. The
paper \cite{f-f-t} solves a long-standing problem by showing that 
not all finite posets are sphere orders. In particular, the poset
$\bm{n}^3$ is not a sphere order for $n$ sufficiently large,
where $\bm{n}$ denotes an $n$-element chain.

For $n=1$ the situation is much simpler. The timelike posets for $n=1$
(i.e., events in $\rr^{1,1}$) are just the posets of dimension 2,
i.e., posets that are an intersection of two chains
\cite[Prop.~2]{meyer}.  Equivalently, they are subposets of $\zz\times
  \zz$ (with the usual product ordering). For characterizations of
  posets of dimension 2, see e.g.\ \cite[{\S}3.3]{trotter}.
 
\section{What permutations of the events are possible?} \label{sec4}
\indent Let $f(n,k)=c(k,k)+c(k,k-1)+\cdots+c(k,k-n)$.  We know from
Theorem~\ref{thm:ssn} that there exist generic spacelike separated
events $\ppk$ in $\rrn$ such that ${\cal O}(\ppk) = f(n,k)$.  Regard
an ordering of these events as a permutation $\pi\in\sk$, the
symmetric group of all permutations of $[k]$. Thus the $k$ events
determine a subset of $\sk$ of cardinality $f(n,k)$. We may further
ask what subsets of $\sk$ of cardinality $f(n,k)$ are possible, and
how many such subsets are there? In general this seems to be a
difficult question, so we will restrict our attention to the case
$n=1$.

Assume then $n=1$, so $\bmv=v\in\rr$. Note that $f(1,k)=1+{k\choose
2}$ by Theorem~\ref{thm:ssn}, since $c(k,k-1)={k\choose 2}$. For the
remainder of this section we continue to assume that
  \beq {\cal O}(\ppk) = 1+{k\choose 2}. \label{eq:op2} \eeq
As $v$ increases from $-1$ to $1$, the order of the 
events will change (as seen from a reference frame moving at velocity
$v$) when $v$ passes through the value
  $$ v_{ij} = \frac{t_i-t_j}{x_i-x_j}, $$
where $\bm{p_i}=(t_i,x_i)$. We thus get a sequence
  $$ \Lambda =\left( \pi_0,\pi_1,\dots,\pi_{{k\choose 2}}\right), $$
of permutations of $1,2,\dots,k$ (in agreement with (\ref{eq:op2})
Theorem~\ref{thm:ssn} in the case $n=1$). Assume without loss of
generality that for $v=0$ the permutation is $12\cdots k$, i.e.,
$t_1<t_2<\cdots<t_k$. Varying the $\bm{p_i}$'s, the sequence $\Lambda$
will change when crossing the surface
  $$ \frac{t_i-t_j}{x_i-x_j} = \frac{t_r-t_s}{x_r-x_s}. $$
Hence the number of different $\Lambda$ is governed by the arrangement  
 \beq  \frac{t_i-t_j}{x_i-x_j} = \frac{t_r-t_s}{x_r-x_s}, 
    \label{eq:hyper} \eeq
$1\leq i<j\leq k$, $1\leq r<s\leq k$, of quadric hypersurfaces in
$\rr^n\times\rr^n$. In particular, the number of regions of this
arrangement is an upper bound for the number of $\Lambda$. We can't be
sure that equality holds since we could have
two different regions lying on the same side of all the hypersurfaces.
Note, however, that if we fix the times $t_1,\dots,t_k$ (or the points
$x_1,\dots,x_k$), then (\ref{eq:hyper}) defines a hyperplane
arrangement ${\cal D}={\cal D}(t_1,\dots,t_k)$. Thus in this situation
the number of different $\Lambda$ is just $r({\cal D})$.  In general
it seems difficult to compute this number. A special case was
considered in a different context in \cite{kott}; see the paragraph
after Theorem~\ref{thm:gt} below.

%

\medskip
\begin{example}
Let $(\bm{p_1},\bm{p_2},\bm{p_3},\bm{p_4}) =
((0,1),(1,6),(2,4),(3,11))$. Then for instance 
  $$ v_{12} =\frac{1-0}{6-1}=\frac 15, $$
and we obtain
  $$ v_{23}<0<v_{34}<v_{12}<v_{14}<v_{24}<v_{13}. $$
Hence $\Lambda=(1324,1234,1243,2143,2413,4213,4231)$.
\end{example}

If $\Lambda=\left( \pi_0,\pi_1,\dots,\pi_{{k\choose 2}}\right)$,
then $\pi_{i+1}$ differs from $\pi_i$ by an adjacent
transposition. Hence 
  $$ \Lambda  =  \pi_0\cdot\left(
  \rho_0,\rho_1,\dots,\rho_{{k\choose 2}}\right) = 
    \left( \pi_0\rho_0,\pi_0\rho_1,\dots,
       \pi_0\rho_{{k\choose 2}}\right), $$
where some $\rho_i=\pi_0^{-1}$, and $\left(
\rho_0,\rho_1,\dots,\rho_{{k\choose 2}}\right)$ is a maximal chain in
the \emph{weak (Bruhat) order} of $\sk$ (see e.g.\
\cite{e-g}\cite{garsia}\cite{rs:rd}). This means that $\rho_0=123\cdots
k$ (the identity permutation), $\rho_{{k\choose 2}}=k\cdots 21$ (the
permutation $w_0\in\sk$ of longest length, i.e., with the most number
  of pairs out of order), and for all $1\leq i\leq
{k\choose 2}$ we have $\rho_i =s_{a_i}\rho_{i-1}$
for some adjacent transposition $s_{a_i}=(a_i,a_i+1)$. It is
well-known (see the previous three references) that the number of
maximal chains in the weak  order of $\sk$ is equal to the
number of standard Young tableaux $f^{(k-1,k-2,\dots,1)}$ of shape
$(k-1,k-2,\dots, 1)$, given by
  $$ f^{(k-1,k-2,\dots,1)} = \frac{{k\choose 2}!}{1^{k-1}\,3^{k-2}\,
      5^{k-3}\cdots (2k-3)^1}. $$
Since $\pi_0=\rho^{-1}$ for some $i$, this gives an upper bound of
    $$ \left( 1+{k\choose 2}\right) f^{(k-1,k-2,\dots,1)} $$ 
for the number of possible $\Lambda(\ppk)$ . Note that if the chain
$\left(\rho_0,\dots,\rho_{{k\choose 2}}\right)$ in the weak order is
achievable, then so is $\sigma\cdot \left(\rho_0,\dots,\rho_{{k\choose
      2}}\right)$ whenever $\sigma=\rho_i^{-1}$ for some $i$, since
$\sigma$ simply specifies which reference frame (or velocity $v$) we
regard as the rest frame ($v=0$). Thus for the problem of
characterizing the possible sequences
$\Lambda(\ppk)=\left(\pi_0,\dots,\pi_{{k\choose 2}}\right)$, we may
assume that $\pi_0=12\cdots k$ or equivalently, $\Lambda(\ppk)$ is a
maximal chain in the weak order of $\sk$.

When $k=3$, it is easy to find examples of all eight sequences
$\sigma\cdot (\rho_0,\rho_1,\rho_2,\rho_3)$ (or of the two maximal chains
in the weak  order of $\mathfrak{S}_3$). To be concrete, these
sequences are
  $$ \begin{array}{c} (123,132,312,321),\ (213,123,132,312),\
  (231,213,123,132)\\
   (321,231,213,123),\ (321,312,132,123),\ (312,132,123,213\\
    (132,123,213,231),\ (123,213,231,321). \end{array} $$
In particular, none of these sequences contain all three of 123, 231,
312, thereby justifying the assertion made at the end of
Section~\ref{sec1}. 

When $k=4$ it can also be checked that all $16$ maximal chains in the
weak order of $\mathfrak{S}_4$ are achievable.  However, for $k=5$ not
all maximal chains occur (see equation (\ref{eq:notal})).  This can be
seen by rephrasing the question of characterizing $\Lambda(\ppk)$ in
terms of earlier work of Goodman and Pollack.  Regard the events
$\bm{p_i}=(t_i,x_i)$ as vectors in $\rr^2$. The order of the events
will be the order they appear when orthogonally projected to a line
$x=C$. We regard this line as having slope $0$. The velocity
$v_{ij}=(t_i-t_j)/(x_i-x_j)$ at 
which the order of the events $\bm{p_i}$ and $\bm{p_j}$ changes is just
the reciprocal of the slope of the line through $\bm{p_i}$ and
$\bm{p_j}$.  It follows that the order of the events in the reference
frame moving at velocity $v$ is just the order in which they appear
when projected to a line of slope $-v$. In other words, as we rotate
the line $t=-x$ (slope $-1$) counterclockwise through an angle
of $90^\circ$ (so it becomes the line $t=x$ of slope 1), the order of
the projections of $\ppk$ on this line 
will change each time the line becomes perpendicular to a line through
two of the points. Regard an ordering of the points $\ppk$ as a
permutation $\pi\in\sk$. It follows that the sequences
$\left(\pi_0,\dots,\pi_{{k\choose 2}}\right)$ of permutations
$\pi_i\in\sk$ obtained in this way from planar configurations of
points are exactly the sequences $\Lambda(\ppk)$. Such sequences were
considered by Goodman and Pollack \cite{g-p} in connection with some
problems of discrete geometry. (They considered lines of all slopes
$\sigma$, not just $-1<\sigma<1$, but this does not produce any
greater generality because we can replace $(t,x)$ with $(t,ax)$ for
$a\gg 0$.) They showed (Theorem~3.3) that all 
maximal chains in the weak order of $\sk$ can occur for $k\leq 4$,
but that for $k=5$ the sequence $(\pi_0,\dots,\pi_{10})$ is not
achievable, where $\pi_i=s_{a_i}\pi_{i-1}$ and
  \beq (a_1,\dots,a_{10}) = (1,3,4,2,1,3,4,2,1,3). \label{eq:notal}
  \eeq 
Hence the same is true for the sequences $\Lambda(\ppk)$.

\section{A classical analogue.} \label{sec5}
\indent A result of Good and Tideman \cite{g-t} may be regarded as a
special case of Theorem~\ref{thm:ssn}. We state their result in a form
involving classical physics so that it is more analogous to
Theorem~\ref{thm:ssn}, though it really has nothing to do with
physics. Suppose that $\ppk\in\rr^n$ (Euclidean space). At time $t=0$
each point $\bm{p_i}$ emits a flash of light. In how many orders can
these events be observed from different points $\bm{x}\in\rr^n$? First
note the fundamental difference between this question and the
situation of Theorem~\ref{thm:ssn}, viz., now we are concerned with
the order in which events are \emph{observed}, not in which they
\emph{occur}. (Of course in classical physics, the order in which
events occur is the same in all reference frames.)

The events $\bmp$ and $\bmq$ are observed simultaneously at points
$\bmx$ on the perpendicular bisector of $\bmp$ and $\bmq$, with
equation
  $$ (\bmp-\bmq)\cdot \bmx=\frac 12\left(|\bmp|^2-|\bmq|^2\right). $$
Hence in analogy with Theorem~\ref{main} we obtain the following
result.

\medskip
\begin{theorem}
\emph{The number of different orders in which $\ppk$ can be observed
  is the number $r({\cal C})$ of regions of the arrangement ${\cal C}$
  with hyperplanes}
  \beq  (\bm{p_i}-\bm{p_j})\cdot \bm{x} = \frac 12\left(|\bm{p_i}|^2-
     |\bm{p_j}|^2\right),\ \ 1\leq i<j\leq n.  \label{eq:pb} \eeq
\end{theorem}
\medskip

This arrangement (\ref{eq:pb})  is a special case of equation
(\ref{eq:hyp}). Moreover, the genericity of $\ppk$ in (\ref{eq:pb}) is
sufficient for genericity in the sense of Theorem~\ref{thm:ssn}. We
therefore obtain the next result.

\medskip
\begin{theorem} \label{thm:gt}
\emph{Let $\ppk$ be generic events in $\rr^n$ occuring at $t=0$. Then
  the number of different orders in which these events can be observed
  at points $\bmx\in\rr^n$ is given by}
  $$ r({\cal C}) = c(k,k)+c(k,k-1)+\cdots+c(k,k-n). $$
\end{theorem}
 
Theorem~\ref{thm:gt} may be restated as determining the number of
regions into which $\rr^n$ is divided by the perpendicular bisectors
of $k$ generic points. This problem was first considered by Good and
Tideman \cite{g-t} in connection with voting theory. They obtained our
Theorem~\ref{thm:gt} by a rather complicated induction argument.
Zaslavsky \cite{zas2} corrected an oversight in the proof of Good and
Tideman and reproved their result by using standard
techniques from the theory of arrangements. Zaslavsky's proof is more
complicated than ours, but he works in a more general context.
Recently Kamiya, Orlik, Takemura, and Terao \cite{kott} considered
additional aspects of Theorem~\ref{thm:gt} in an analysis of ranking
patterns, in particular, enumerating the number of sets of orders that
can occur by varying the points $\ppk$. 

\medskip
\textsc{Acknowledgement.} I am grateful to a person from the National
Security Agency whose name I cannot recall for suggesting to me the
topic of this paper and for mentioning that it is connected with
arrangments of hyperplanes. I am also grateful to Daniel Freedman for
some helpful discussions on special relativity.

\pagebreak

\end{document}